\documentclass[11pt,reqno]{amsart}

\usepackage{amssymb,amsxtra}
\usepackage{concrete}
\usepackage{eucal}
\usepackage[T1]{fontenc}
\usepackage{textcomp}
\usepackage{latexsym}
\usepackage{epic,eepic}

\usepackage[cmtip,matrix,arrow,curve]{xy}
\CompileMatrices
\newdir{((}{{\lhook\kern-.1em}}
\newdir{))}{{\rhook\kern-.1em}}
\newdir{ >}{{}*!/-5pt/@{>}}

\newcommand\note[1]%
{$^\dagger$\marginpar{\footnotesize{$^\dagger${#1}}}}

\renewcommand\dots{\relax\ifmmode\ldots\else$\,\ldots\,$\fi}

\divide\count253 by 60 %hours after midnight, rounded down
\divide\count255 by 60
\multiply\count255 by 60 %minutes after midnight, rounded down
\advance\count254 by -\count255 %minutes after the hour

\def\today{\number\day\space\ifcase\month\or January\or February\or
March\or April\or May\or June\or July\or August\or September\or
October\or November\or December\fi\space\number\year}

\def\hour{\ifnum\count253<10
0\number\count253\else\number\count253\fi}

\def\minute{\ifnum\count254<10
0\number\count254\else\number\count254\fi}

\swapnumbers
\newtheorem{theorem}{Theorem}[section]
\newtheorem*{theorem*}{Theorem}

\newtheorem{proposition}[theorem]{Proposition}
\newtheorem{corollary}[theorem]{Corollary}

\theoremstyle{definition}

\newtheorem{remark}[theorem]{Remark}
\newtheorem{remarks}[theorem]{Remarks}

\newcommand\lie{\mathfrak}

\renewcommand{\t}{\mathfrak{t}}
\newcommand{\g}{\mathfrak{g}}
\newcommand{\h}{\mathfrak{h}}
\newcommand{\n}{\mathfrak{n}}
\renewcommand{\b}{\mathfrak{b}}
\newcommand{\p}{\mathfrak{p}}

\newcommand\bb[1]{{\text{\bf#1}}}
\newcommand\N{\bb{N}}
\newcommand\Z{\bb{Z}} 
\newcommand\Q{\bb{Q}}
\newcommand\R{\bb{R}} 
\newcommand\C{\bb{C}}

\newcommand\CP{\bb{CP}}

\newcommand\ca{\mathcal}

\newcommand\func[1]{\operatorname{\mathrm{#1}}}

\newcommand\Ad{\func{Ad}}

\newcommand\Hom{\func{Hom}}

\newcommand\Int{\func{Int}}

\newcommand\Proj{\func{Proj}}
\newcommand\pr{\func{pr}}

\renewcommand\Re{\func{Re}}
\renewcommand\Im{\func{Im}}

\newcommand\group[1]{{\text{\bf#1}}}

\newcommand\U{\group{U}}

\newcommand\norm[1]{\lVert#1\rVert}

\newcommand\inner[1]{\langle#1\rangle}

\newcommand\quot[1][\kern.3ex]{/\kern-.7ex/_{\kern-.4ex#1}}
\newcommand\bigquot[1][\,\,]{\big/\kern-.85ex\big/_{\!\!#1}}

\newcommand\powl{[\kern-.3ex[}
\newcommand\powr{]\kern-.3ex]}
\newcommand\bigpowl{\bigl[\kern-.6ex\bigl[}
\newcommand\bigpowr{\bigr]\kern-.6ex\bigr]}

\newcommand\restrict{\mathop{|}}

\newcommand\sur{\mathrel{\to\kern-1.8ex\to}}
\newcommand\iso{\mathrel{\hookrightarrow\kern-1.8ex\to}}

\newcommand\longto{\longrightarrow}
\newcommand\longhookrightarrow{\lhook\joinrel\longrightarrow}

\newcommand\longinj{\longhookrightarrow}
\newcommand\longsur{\mathrel{\longrightarrow\kern-1.8ex\to}}
\newcommand\longiso{\mathrel{\longhookrightarrow\kern-1.8ex\to}}

\newcommand\sq{\sqrt{-1}}

\renewcommand\subset{\subseteq}
\renewcommand\supset{\supseteq}

\newcommand\impl{_{\mathrm{impl}}}
\renewcommand\ss{^{\mathrm{ss}}} 

\begin{document}

%%%%%%%%%%%%%%%%%%%%%%%%%%%%%%%%%%%%%%%%%%%%%%%%%%%%%%%%%%%%%%%%%%%%%%%%
%%%%%%%%%%%%%%%%%%%%%%%%%%%%%%%%%%%%%%%%%%%%%%%%%%%%%%%%%%%%%%%%%%%%%%%%

\title[Convexity for $B$-varieties]{Convexity theorems for varieties
invariant under a Borel subgroup}

\dedicatory{To Bob MacPherson, our teacher and friend}

\author{Victor Guillemin}

\address{Department of Mathematics, Massachusetts Institute of
Technology, Cambridge, Massachusetts 02139-4307}

\email{vwg@math.mit.edu}

\author{Reyer Sjamaar}

\address{Department of Mathematics, Cornell University, Ithaca, New
York 14853-4201}

\email{sjamaar@math.cornell.edu}

\date{25 January 2006}


\begin{abstract}
Atiyah proved that the moment map image of the closure of an orbit of
a complex torus action is convex.  Brion generalized this result to
actions of a complex reductive group.  We extend their results to
actions of a maximal solvable subgroup.
\end{abstract}

\maketitle

\tableofcontents

%%%%%%%%%%%%%%%%%%%%%%%%%%%%%%%%%%%%%%%%%%%%%%%%%%%%%%%%%%%%%%%%%%%%%%%%
%%%%%%%%%%%%%%%%%%%%%%%%%%%%%%%%%%%%%%%%%%%%%%%%%%%%%%%%%%%%%%%%%%%%%%%%

%%%%%%%%%%%%%%%%%%%%%%%%%%%%%%%%%%%%%%%%%%%%%%%%%%%%%%%%%%%%%%%%%%%%%%%%
\section*{Introduction}
%%%%%%%%%%%%%%%%%%%%%%%%%%%%%%%%%%%%%%%%%%%%%%%%%%%%%%%%%%%%%%%%%%%%%%%%

Let $M$ be a compact K\"ahler manifold and $L$ a Hermitian line bundle
over $M$ whose curvature form is the K\"ahler form.  In addition, let
$G$ be a compact Lie group and $\tau$ a K\"ahlerian action of $G$ on
$(M,L)$.  From this action one gets a moment map $\Phi\colon M\to\g^*$
and a holomorphic action $\tau_\C$ of the complexified Lie group
$G_\C$ on $(M,L)$.  In particular, any Borel subgroup $B$ of $G_\C$
acts on $M$ and $L$, and the main goal of this paper is to prove a
Kirwan type convexity theorem for the moment images of $B$-invariant
subvarieties of $M$.  Let $T$ be the maximal torus of $G$ contained in
$B$ and let $\t^*_+$ be the closed Weyl chamber in $\t^*$ which is
positive with respect to $B$.  Our main result, which is described in
more detail in Section \ref{section;borel}, asserts

\begin{theorem*}
If $X$ is a $B$-invariant irreducible closed analytic subvariety of
$M$, the moment image of $X$ in $\g^*$ intersects $\t^*_+$ in a convex
polytope.
\end{theorem*}

This result is a $B$-analogue of well-known convexity theorems of
Atiyah and Brion, which we will review in Section
\ref{section;atiyah-brion}, partly as a motivation for our result and
partly because the techniques used by Brion to prove his
$G_\C$-convexity theorem will be used in our proof as well.

One other major ingredient in our proof is a $B$-analogue of a
fundamental result in geometric invariant theory.  This $B$-analogue
asserts that if $s$ is a global holomorphic section of $L$ which is an
eigensection for the action of $B$ and which does not vanish at a
point $x$, then the norm of $s$, restricted to the closure of the
$B$-orbit through $x$, takes its maximum at a point $y$ such that
$\Phi(y)$ is in the chamber $\t^*_+$.  This theorem and variants of it
will be discussed in Section \ref{section;norm}.  The final ingredient
of our proof, a rationality statement concerning the moment image, is
explained in Section \ref{section;dense}.

Our main result has some interesting applications to ``complex Morse
theory''.  It is well known that a generic component of the moment map
is a Morse-Bott function, and Carrell and Sommese have shown that the
unstable manifolds of this Morse-Bott function are complex
submanifolds and that their closures are analytic subvarieties.
Moreover, as we show in Section \ref{section;schubert}, for a
``dominant'' choice of a moment map component these varieties are
$B$-invariant.  (For instance, for a flag variety the unstable
manifolds are the Bruhat cells: the $B$-orbits.)  Hence their moment
images intersected with $\t^*_+$ are convex polytopes.  In other
words, the Morse-Bruhat-Carrell-Sommese decomposition of $M$ gets
reflected in a ``stratification'' of the Kirwan polytope.

%%%%%%%%%%%%%%%%%%%%%%%%%%%%%%%%%%%%%%%%%%%%%%%%%%%%%%%%%%%%%%%%%%%%%%%%
\section{Convexity results of Atiyah and Brion}
\label{section;atiyah-brion}
%%%%%%%%%%%%%%%%%%%%%%%%%%%%%%%%%%%%%%%%%%%%%%%%%%%%%%%%%%%%%%%%%%%%%%%%

The main result of this paper is a nonabelian generalization of a
theorem about torus actions proved by Atiyah in
\cite{atiyah;convexity-commuting}.  Let $(M,\omega)$ be a compact
K\"ahler manifold and $\tau$ a K\"ahlerian action of the torus
$T=(S^1)^n$ on $M$.  From $\tau$ one gets a holomorphic action
$\tau_\C$ of the complex torus $T_\C=(\C^*)^n$.  If $\tau$ has fixed
points, there is a moment map $\Phi\colon M\to\t^*$, and Atiyah's
theorem asserts:

\begin{theorem}\label{theorem;atiyah}
Let $X=\overline{T_\C x}$ be the closure of the $T_\C$-orbit through a
point $x\in M$.  Then the moment image of $X$ is a convex polytope.
More explicitly, $\Phi(X)$ is the convex hull of the set $\Phi(X^T)$,
where $X^T$ denotes the set of fixed points of $\tau\restrict X$.
\end{theorem}

If $\omega$ is the curvature form of a $T$-equivariant Hermitian
holomorphic line bundle $L$ on $M$, there is an alternative
description of $\Phi(X)$.  The space of global holomorphic sections
$\Gamma(M,L)$ decomposes under the action of $T$ into a direct sum of
weight spaces,
$\Gamma(M,L)=\bigoplus_{\lambda\in\Lambda}\Gamma(M,L)_\lambda$, where
$\Lambda=\Hom(T,\U(1))$ is the weight lattice of $T$.  (As usual, we
will identify $\Lambda$ with a subgroup of $\t^*$ by identifying the
weight $\lambda$ with the functional $T_1\lambda/(2\pi\sq)$.)
Atiyah's description of $\Phi(X)$ is equivalent to:

\begin{theorem}\label{theorem;atiyah-weight}
Let $x\in M$ and let $X=\overline{T_\C x}$.  Suppose $L$ is generated
by its global sections.  Then $\Phi(X)$ is the convex hull of the
subset
$$
\{\,\lambda\in\Lambda\mid\exists
s\in\Gamma(M,L)_\lambda,\;s(x)\ne0\,\}
$$
of $\t^*$.
\end{theorem}

It is this version of Atiyah's theorem which we'll generalize below.
First, however, we mention some of its consequences.  

\begin{corollary}\label{corollary;hull}
Let $X$ be an arbitrary $T$-invariant irreducible closed subvariety of
$M$ and let $i\colon X\to M$ be the inclusion map.
\begin{enumerate}
\item\label{item;hull}
$\Phi(X)$ is the convex hull in $\t^*$ of the set
$$
\{\,\lambda\in\Lambda\mid\exists s\in\Gamma(M,L)_\lambda,\;i^*s\ne0\,\}.
$$
\item\label{item;fix-hull}
$\Phi(X)$ is the convex hull in $\t^*$ of the set $\Phi(X^T)$.
\item\label{item;integral}
$\Phi(X)$ is an integral convex polytope.
\item
The set consisting of all $x$ in $X$ for which $\Phi(\overline{T_\C
x})=\Phi(X)$ is nonempty and Zariski open in $X$.
\item\label{item;open}
The collection of polytopes $\Phi(X)$, where $X$ ranges over all
$T_\C$-invariant irreducible closed subvarieties of $M$, is finite.
\end{enumerate}
\end{corollary}

A special case of Corollary \ref{corollary;hull} is a convexity
theorem of Atiyah for Schubert varieties.  Suppose $T$ is a Cartan
subgroup of a connected compact Lie group, $G$.  Let $G_\C$ be the
complexification of $G$, let $B\supset T$ be a Borel subgroup of
$G_\C$ and let $M$ be the full flag variety $G_\C/B$.  For any
$\mu\in\Lambda$ let $\C_\mu$ be the one-dimensional $B$-module with
weight $\mu$ and let $L_\mu=G_\C\times^B\C_\mu$ be the corresponding
$G_\C$-homogeneous line bundle over $M$.  The curvature form
$\omega_\mu$ (with respect to the $G$-invariant Hermitian metric on
$L_\mu$ defined by the absolute value on $\C$) is nondegenerate if and
only if the weight $\mu$ is regular.  It is well known that
$\omega_\mu$ is K\"ahler if and only if $\mu\in-\Int\t^*_+$, where
$\t^*_+$ is the closed Weyl chamber in $\t^*$ which is positive with
respect to $B$.  A moment map $\Phi_\mu$ for the $G$-action is found
by composing the maps
\begin{equation}\label{equation;flag}
M=G_\C/B\longto G/T\longto\g^*,
\end{equation}
where the first map is the inverse of the diffeomorphism $G/T\to M$
induced by the inclusion $G\to G_\C$ and the second map is defined by
$gT\mapsto g\mu$.  In fact, $\Phi_\mu$ is a $G$-equivariant
symplectomorphism from $M$ onto the coadjoint orbit through $\mu$.
Let $W=N_G(T)/T$ be the Weyl group of $(G,T)$.  The $T$-fixed points
in $M$ are precisely the cosets $wB$ with $w\in W$.  Every orbit of
the left $B$-action on $M$ passes through a unique $T$-fixed point.
The \emph{Schubert variety} $X_w$ is by definition the closure of the
$B$-orbit through the $T$-fixed point $wB$.  The \emph{Bruhat order}
on $W$ is defined by $v\le w$ if and only if $X_v\subset X_w$.  Thus
$(X_w)^T=\{\,vB\mid v\le w\,\}$.  Since the $T$-moment map $M\to\t^*$
is just the composition of $\Phi_\mu$ and the projection
$\g^*\to\t^*$, one obtains from Corollary
\ref{corollary;hull}\eqref{item;fix-hull}:

\begin{theorem}[\cite{atiyah;convexity-commuting}]
\label{theorem;schubert}
The projection onto $\t^*$ of the Schubert variety $X_w$ is equal to
the convex hull of the set $\{\,v\mu\mid v\le w\,\}$.
\end{theorem}

Theorem \ref{theorem;atiyah-weight} and its various corollaries have
been generalized by Brion to nonabelian groups in the setting of
complex projective varieties.  Namely let $G$, $T$ and $B$ be as above
and let $N=[B,B]$ be the unipotent radical of $B$.  Let $\tau$ be an
action of $G$ on $\CP^n$ coming from a unitary representation of $G$
on $\C^{n+1}$.  Then $G$ acts on the polynomial algebra
$S=\bigoplus_{r\ge0}S_r$, where $S_r=\Gamma(\CP^n,L^r)$ and
$L=\ca{O}(1)$ is the canonical hyperplane bundle.  Let
$\Lambda_+=\Lambda\cap\t^*_+$ be the semigroup of dominant weights of
$(G,T)$.  For $\lambda\in\Lambda_+$ and $r\in\N$ let $S_{\lambda,r}$
be the isotypical $G$-submodule of $S_r$ of highest weight $\lambda$.
The set $S^N$ consisting of all $N$-invariant polynomials is a
subalgebra of $S$ which is preserved by the action of $T$.  Let us put
$S^N_{\lambda,r}=(S_{\lambda,r})^N$; then
$$
S^N=\bigoplus\bigl\{\,S^N_{\lambda,r}\bigm|
(\lambda,r)\in\Lambda_+\times\N\,\bigr\}
$$
is a bigrading of the algebra $S^N$ by weight and degree.  Let $X$ be
any $G_\C$-invariant irreducible closed subvariety $X$ of $\CP^n$ and
let $I(X)$ be the homogeneous ideal of $S$ consisting of all
polynomials vanishing on $X$.  Then $I(X)^N$ is a homogeneous ideal of
$S^N$ and the quotient algebra $A(X)=S^N/I(X)^N$ is
$\Lambda_+\times\N$-graded.  Brion defines
$$
\ca{C}(X)=\{\,\lambda\in\Lambda\otimes\Q
\mid\text{$A(X)_{r\lambda,r}\ne0$ for some $r>0$}\,\}.
$$
In other words, a rational weight $\lambda\in\Lambda\otimes\Q$ is in
$\ca{C}(X)$ if and only if there exist a positive integer $r$ such
that $r\lambda$ is integral and a section $s\in S^N_{r\lambda,r}$ such
that $s$ does not vanish on $X$.  The set $\ca{C}(X)$ is contained in
the intersection of the chamber $\t^*_+$ with the $\Q$-vector space
$\Lambda\otimes\Q$.  The assumption that $X$ is irreducible implies
that $A(X)$ has no zero divisors and hence that $\ca{C}(X)$ is convex
over the rationals.  A theorem of Hadziev and Grosshans (see e.g.\
\cite[Kapitel III.3.2]{kraft;geometrische} says that the algebra $S^N$
is finitely generated.  Since $A(X)$ is the image of the map
$S^N\to(S/I(X))^N$ induced by the quotient map $S\to S/I(X)$, $A(X)$
is likewise finitely generated.  By selecting a finite system of
homogeneous generators $\{\,a_i\in
A(X)_{\lambda_i,r_i}\mid\text{$i=1$, $2$,\dots, $k$}\,\}$, one sees
immediately that $\ca{C}(X)$ is the convex hull of the rational
weights $\lambda_1/r_1$, $\lambda_2/r_2$,\dots, $\lambda_k/r_k$.  Thus
$\ca{C}(X)$ is a convex polytope over $\Q$.  Now let
$\Phi\colon\CP^n\to\g^*$ be the moment map for the action of $G$ on
$\CP^n$.  Brion proves:

\begin{theorem}[\cite{brion;image}]\label{theorem;brion}
Let $X$ be a $G_\C$-invariant irreducible closed subvariety of
$\CP^n$.  Then the intersection $\Delta(X)=\Phi(X)\cap\t^*_+$ is equal
to the closure of $\ca{C}(X)$, and $\ca{C}(X)$ is equal to the set of
rational points in $\Delta(X)$.  Hence $\Delta(X)$ is a rational
convex polytope.
\end{theorem}

(This result builds on prior work by Ness and Mumford
\cite{ness;stratification}; cf.\ also Guillemin and Sternberg
\cite{guillemin-sternberg;convexity;;1982}.)  In particular, for every
$x\in\CP^n$ the set $\Delta(\overline{G_\C x})$ is a rational convex
polytope.  Brion also proves nonabelian analogues of the other
assertions in Corollary \ref{corollary;hull}.  (However, in contrast
to the abelian case, $\Delta$ is not necessarily an \emph{integral}
polytope.  This is why the powers of $L$ must be incorporated in the
definition of the set $\ca{C}(X)$.)
% Observe also that Theorem \ref{theorem;atiyah-weight} is strictly
% speaking not a special case of Theorem \ref{theorem;brion}.
% However, Theorem \ref{theorem;atiyah-weight} can be easily deduced
% from Theorem \ref{theorem;brion} together with Theorem
% \ref{theorem;atiyah}.

%%%%%%%%%%%%%%%%%%%%%%%%%%%%%%%%%%%%%%%%%%%%%%%%%%%%%%%%%%%%%%%%%%%%%%%%
\section{Convexity theorems for $B$-varieties}\label{section;borel}
%%%%%%%%%%%%%%%%%%%%%%%%%%%%%%%%%%%%%%%%%%%%%%%%%%%%%%%%%%%%%%%%%%%%%%%%

Our main observation is that these results of Brion are also true for
$B$-invariant subvarieties of $G_\C$-manifolds.  Namely let $M$ be a
compact complex manifold, $L \to M$ a positive Hermitian holomorphic
line bundle with Hermitian connection $\nabla$ and curvature form
$\omega$, and $\tau$ an action of $G$ on $L$ by line bundle
automorphisms which preserve the complex structure on $M$ and the
Hermitian structure on $L$.  The group of \emph{all} line bundle
automorphisms which preserve the holomorphic structure on $L$ is a
complex Lie group; and therefore the $G$-action extends uniquely to a
$G_\C$-action by holomorphic line bundle automorphisms on $L$.
Moreover, since the action of $G$ on $L$ preserves its Hermitian
structure, the action of $G$ on $M$ preserves $\omega$.  The space
$C^\infty(M,L)$ of all smooth sections of $L$ is a $G$-module in a
natural way.  Let $\ca{L}(\xi)s$ denote the Lie derivative of a
section $s\in C^\infty(M,L)$ along $\xi\in\g$.  As observed by Kostant
\cite[Theorem 4.3.1]{kostant;quantization-prequantization}, the
first-order operators $\nabla(\xi_M)$ and $\ca{L}(\xi)$ have the same
principal symbol, and the zeroth order operator
$\ca{L}(\xi)-\nabla(\xi_M)$ is given by multiplication by an
imaginary-valued function.  The map $\phi\colon\g\to C^\infty(M,\R)$
defined by
\begin{equation}\label{equation;nabla}
\phi(\xi)=\frac1{2\pi\sq}(\ca{L}(\xi)-\nabla(\xi_M))
\end{equation}
satisfies $d\phi(\xi)=\iota(\xi_M)\omega$ and is $G$-equivariant.
Therefore the map $\Phi\colon M\to\g^*$ defined by
$\inner{\Phi,\xi}=\phi(\xi)$ is a moment map for the $G$-action on
$M$.

As in Section \ref{section;atiyah-brion} we consider the algebra of
sections $S=\bigoplus_{r\ge0}\Gamma(M,L^r)$ and its subalgebra $S^N$,
which has a grading
$$
S^N=\bigoplus\bigl\{\,S^N_{\lambda,r}\bigm|
(\lambda,r)\in\Lambda_+\times\N\,\bigr\}.
$$
We introduce, for any $B$-invariant irreducible closed analytic
subvariety $X$ of $M$, the homogeneous ideal $I(X)$ of all sections
vanishing on $X$, the quotient algebra $A(X)=S^N/I(X)^N$, and the set
of rational weights
$$
\ca{C}(X)=\{\,\lambda\in\Lambda\otimes\Q
\mid\text{$A(X)_{r\lambda,r}\ne0$ for some $r>0$}\,\}.
$$
(Sometimes we shall write $\ca{C}(X,L)$ instead of $\ca{C}(X)$ to
emphasize the dependence of this set on the linearization $L$.)  As
before, $\ca{C}(X)$ is contained in the chamber $\t^*_+$ and is the
convex hull of a finite subset of $\Lambda\otimes\Q$.  Now define
$\Delta(X)=\Phi(X)\cap\t^*_+$.  In Sections \ref{section;norm} and
\ref{section;dense} we will prove the following assertion, which is
the main result of this paper.

\begin{theorem}\label{theorem;borel}
Let $X$ be a $B$-invariant irreducible closed subvariety of $M$.  Then
$\Delta(X)$ is equal to the closure of $\ca{C}(X)$, and $\ca{C}(X)$ is
equal to the set of rational points in $\Delta(X)$.  Hence $\Delta(X)$
is a rational convex polytope.
\end{theorem}

\begin{remarks}\label{remarks}
\begin{list}{(\roman{enumi})}{\leftmargin0pt\labelsep1ex%
\labelwidth-1ex\itemindent1em\usecounter{enumi}}
\item
By Kodaira's embedding theorem, $M$ embeds holomorphically into the
projective space associated to $\Gamma(M,L^r)^*$ for sufficiently
large $r$.  It is easy to see that the Kodaira map is
$G_\C$-equivariant.  Therefore $M$ is a projective $G_\C$-variety.
However, the Kodaira map is rarely symplectic.  Thus the symplectic
form $\omega$ is usually not the restriction of the Fubini-Study form,
nor is the moment map $\Phi$ usually the restriction of the standard
moment map on projective space.
\item
Here is a heuristic argument in favour of Theorem \ref{theorem;borel}.
Let $M\impl$ be the imploded cross-section of $M$.  This is a certain
stratified symplectic space with a Hamiltonian $T$-action on the
strata.  (See \cite{guillemin-jeffrey-sjamaar}.)  It can be identified
with $\Proj S^N$, the ``quotient'' of $M$ by the $N$-action.  Any
$B$-invariant subvariety $X$ of $M$ maps to a $T_\C$-invariant
subvariety $Y$ of $M\impl$, whose $T$-moment map image is equal to
$\Delta(X)$.  But the image of $Y$ is convex by Atiyah's theorem.
Therefore $\Delta(X)$ is convex.  It is possible to give a rigorous
argument along these lines, but we will present a shorter and more
direct proof.
\item
In contrast to the $G$-invariant case, the intersection
$\Phi(X)\cap\t^*$ is not Weyl group invariant.  We do not know how
$\Phi(X)$ intersects the chambers of $\t^*$ other than $\t^*_+$.  See
however Remark \ref{remark;chamber} for some partial information.
\item\label{item;ideal}
If $X$ is $G_\C$-invariant, the map $S^N\to(S/I(X))^N$ is surjective
(this follows from the reductivity of $G_\C$) and thus
$A(X)=(S/I(X))^N$.  However, this may fail if $X$ is merely
$B$-invariant, in which case the algebra $(S/I(X))^N$ may have
undesirable properties.
\end{list}
\end{remarks}

Theorem \ref{theorem;borel} implies that $\Delta(\overline{Bx})$ is a
rational convex polytope for every point $x\in M$.  For those points
at which $B$ acts freely, a somewhat different description of these
polytopes can be found in \cite[Section~4.5]%
{guillemin-sjamaar;convexity-properties-hamiltonian}, where we prove

\begin{theorem}\label{theorem;nilpotent}
Let $x\in M$ and suppose that $B$ acts freely at $x$.  Then
\begin{equation}\label{equation;intersect}
\Delta(\overline{Bx})=\t^*_+\cap\bigcap_{n\in N}\Phi_T(\overline{T_\C
nx}),
\end{equation}
where $\Phi_T\colon M\to\t^*$ is the $T$-moment map, i.e.\ the
composition of $\Phi$ with the projection $\g^*\to\t^*$.
\end{theorem}

(In [\emph{loc.\ cit.}, Theorem 4.5.1] $B$ denotes the \emph{negative}
Borel relative to the chamber $\t^*_+$.  The difference arises from
the opposite sign convention for the moment map adopted there.)  A
similar result, for moment polytopes of $G_\C$-invariant subvarieties,
was proved by Franz \cite{franz;moment}.  Notice that by Atiyah's
theorem each of the sets $\Phi_T(\overline{T_\C nx})$ is a convex
polytope and only a finite number of distinct polytopes can occur in
the intersection \eqref{equation;intersect}.  Hence
\eqref{equation;intersect} is a convex polytope.  However, our proof
of Theorem \ref{theorem;borel} will be independent of Theorem
\ref{theorem;nilpotent}.

The proof of Theorem \ref{theorem;borel} involves the following
result, which is of interest in its own right:

\begin{theorem}\label{theorem;maximum}
Let $X$ be a $B$-invariant irreducible closed subvariety of $M$ and
let $i\colon X\to M$ be the inclusion map.  Suppose $s\in
S^N_{r\lambda,r}$ satisfies $i^*s\ne0$.  Let $x$ be a point on $X$
where $\norm{i^*s}$ takes its maximum value.  Then $\Phi(x)=\lambda$.
\end{theorem}

We mention a few consequences of Theorem \ref{theorem;borel}.  As in
the torus and reductive cases, $\Delta(X)=\Delta(\overline{Bx})$ for
``most'' $x$ in $X$:

\begin{corollary}\label{corollary;orbit}
Let $X$ be a $B$-invariant irreducible closed subvariety of $M$.  The
set $U$ consisting of all $x$ in $X$ for which
$\Delta(X)=\Delta(\overline{Bx})$ is nonempty and Zariski open in $X$.
\end{corollary}

\begin{proof}
It follows from Theorem \ref{theorem;borel} that
\begin{equation}\label{equation;open}
U=\{\,x\in X\mid\ca{C}(X)=\ca{C}(\overline{Bx})\,\}.
\end{equation}
Let $\lambda_1$, $\lambda_2$,\dots, $\lambda_l$ be the vertices of
$\ca{C}(X)$.  For each $i$ put
\begin{multline*}
U_i=\{\,x\in X\mid\lambda_i\in\ca{C}(\overline{Bx})\,\}\\
=\bigl\{\,x\in X\bigm|\text{$s(x)\ne0$ for some $r>0$ and some $s\in
S^N_{r\lambda_i,r}$}\,\bigr\}.
\end{multline*}
It follows from \eqref{equation;open} that $U=U_1\cap
U_2\cap\cdots\cap U_l$.  Each $U_i$ is nonempty and Zariski open, $X$
is irreducible, and therefore $U$ is nonempty and Zariski open.
\end{proof}

By induction on the dimension this implies the following result.

\begin{corollary}\label{corollary;finite}
The collection of polytopes $\Delta(X)$, where $X$ ranges over all
$B$-invariant irreducible closed subvarieties of $M$, is finite.
\end{corollary}

Corollary \ref{corollary;orbit} applies to $G_\C$-invariant
subvarieties $X$.  In particular:

\begin{corollary}\label{corollary;reductive-borel}
The collection of $x\in M$ for which
$\Delta(\overline{Bx})=\Delta(\overline{G_\C x})=\Delta(M)$ is
nonempty and Zariski open.
\end{corollary}

We refer to the following properties as the \emph{lower
semi-continuity} of the function $x\mapsto\Delta(\overline{Bx})$.

\begin{corollary}\label{corollary;semicontinuous}
Let $x$ be an arbitrary point of $M$ and define
\begin{align*}
M_{\ge x}&=\{\,y\in
M\mid\Delta(\overline{By})\supset\Delta(\overline{Bx})\,\},\\
M_{\le x}&=\{\,y\in
M\mid\Delta(\overline{By})\subset\Delta(\overline{Bx})\,\},\\
M_x&=\{\,y\in M\mid\Delta(\overline{By})=\Delta(\overline{Bx})\,\}.
\end{align*}
Then $M_{\ge x}$ is open, $M_{\le x}$ is closed, and $M_x$ is locally
closed in the Zariski topology on $M$.
\end{corollary}

\begin{proof}
Theorem \ref{theorem;borel} shows that $M_{\ge x}=\{\,y\in
M\mid\ca{C}(\overline{By})\supset\ca{C}(\overline{Bx})\,\}$ (and
analogous identities for $M_{\le x}$ and $M_x$).  The proof that
$M_{\ge x}$ is open is now similar to the proof of Corollary
\ref{corollary;orbit}.  (Use the vertices of $\ca{C}(\overline{Bx})$
instead of those of $\ca{C}(X)$.)  The closedness of $M_{\le x}$ then
follows from the observation that
$$
M_{\le x}=M-\bigcup\{\,M_{\ge y}\mid y\not\in M_{\le x}\,\}.
$$
Finally, $M_x=M_{\ge x}\cap M_{\le x}$, so $M_x$ is locally closed.
\end{proof}

%%%%%%%%%%%%%%%%%%%%%%%%%%%%%%%%%%%%%%%%%%%%%%%%%%%%%%%%%%%%%%%%%%%%%%%%
\section{The norm of an $N$-invariant holomorphic section}
\label{section;norm}
%%%%%%%%%%%%%%%%%%%%%%%%%%%%%%%%%%%%%%%%%%%%%%%%%%%%%%%%%%%%%%%%%%%%%%%%

In this section we will prove Theorem \ref{theorem;maximum} and part
of Theorem \ref{theorem;borel}.  The group $G_\C$ acts complex
linearly on the space of smooth global sections $C^\infty(M,L)$.
Therefore $\ca{L}(\sq\,\xi)s=\sq\,\ca{L}(\xi)s$ for all $s\in
C^\infty(M,L)$.  Let $J$ be the complex structure on $M$.  Since $G$
acts holomorphically on $M$, we have $(\sq\,\xi)_M=J\xi_M$ for all
$\xi\in\g$.  Now suppose $s$ is a \emph{holomorphic} section of $L$.
Then, by the Cauchy-Riemann equation,
$\nabla(J\xi_M)s=\sq\,\nabla(\xi_M)s$.  Thus from
\eqref{equation;nabla} we get
\begin{equation}\label{equation;nabla-i}
\ca{L}(\sq\,\xi)s-\nabla\bigl((\sq\,\xi)_M\bigr)s
=\sq\,\bigl(\ca{L}(\xi)-\nabla(\xi_M)\bigr)s=-2\pi\,\phi(\xi)s.
\end{equation}
For arbitrary $\xi\in\g_\C$ let us write $\Re\xi=(\xi+\bar{\xi})/2$
and $\Im\xi=(\xi-\bar{\xi})/2\sq\,$, where $\bar{\xi}$ is the complex
conjugate of $\xi$ relative to the compact real form $\g$ of $\g_\C$.
Also let us write $\phi_\C\colon\g_\C\to C^\infty(M,\C)$ for the
complexification of the moment map $\phi\colon\g\to C^\infty(M,\R)$,
given by $\phi_\C(\xi)=\phi(\Re\xi)+\sq\,\phi(\Im\xi)$ for
$\xi\in\g_\C$.  Then \eqref{equation;nabla} and
\eqref{equation;nabla-i} imply
$$
\ca{L}(\xi)s-\nabla(\xi_M)s=2\pi\sq\,\phi_\C(\xi)s
$$
for all $\xi\in\g_\C$ and all $s\in\Gamma(M,L)$.  Now suppose that,
for some $\xi\in\g_\C$, $s$ is an eigensection for the operator
$\ca{L}(\xi)$ with eigenvalue $c\in\C$.  Then
$$
\nabla(\xi_M)s=\bigl(c-2\pi\sq\,\phi_\C(\xi)\bigr)s,
$$
and hence
$$
\ca{L}(\xi_M)\norm{s}^2
%=\inner{\nabla_{\xi_M}s\mid s}+\inner{s\mid\nabla_{\xi_M}s}
=2\Re\inner{\nabla(\xi_M)s\mid s}=\bigl(2\Re
c+4\pi\,\phi(\Im\xi)\bigr)\norm{s}^2,
$$
where $\inner{{\cdot}\mid{\cdot}}$ denotes the Hermitian inner product
on $L$ and $\norm{\cdot}$ the associated (pointwise) norm.  This
proves the following assertion.

\begin{proposition}\label{proposition;length}
Let $s$ be a global holomorphic section of $L$.  Let $\lie{k}$ be a
complex Lie subalgebra of $\g_\C$ and suppose $s$ transforms under
$\lie{k}$ according to a character $\chi\colon\lie{k}\to\C$.  Then
$$
\ca{L}(\xi_M)\norm{s}^2
=\bigl(2\Re\chi(\xi)+4\pi\,\phi(\Im\xi)\bigr)\norm{s}^2
$$
for all $\xi\in\lie{k}$.
\end{proposition}

\begin{remarks}\label{remark;norm}
\begin{list}{(\roman{enumi})}{\leftmargin0pt\labelsep1ex%
\labelwidth-1ex\itemindent1em\usecounter{enumi}}
\item
The same formula holds for sections of $L^r$, provided we replace
$\phi$ with $r\phi$.
\item\label{item;subvariety}
The formula holds not only for global holomorphic sections, but also
for holomorphic sections defined over a $\lie{k}$-invariant analytic
subvariety of $M$.
\end{list}
\end{remarks}

There are two important special cases of this identity.  First let
$\lie{k}=\g_\C$ and $\chi=0$.  Then we get, for every $G$-invariant
holomorphic section $s$ and $\xi\in\g$,
\begin{equation}\label{equation;semistable}
\ca{L}\bigl((\sq\,\xi)_M\bigr)\norm{s}^2=4\pi\,\phi(\xi)\norm{s}^2.
\end{equation}
Recall that $x\in M$ is \emph{semistable} if there exist $r>0$ and
$s\in\Gamma(M,L^r)^G$ such that $s(x)\ne0$.  Let $x$ be semistable and
$Y=\overline{G_\C x}$.  If $y\in Y$ is a point where $\norm{s\restrict
Y}$ attains its maximum, then \eqref{equation;semistable} implies
$\Phi(y)=0$.  This proves the inclusion ``$\subset$'' of the following
theorem, in which $M\ss$ denotes the set of semistable points.

\begin{theorem}\label{theorem;semistable}
$M\ss=\bigl\{\,x\in M\bigm|0\in\Phi(\overline{G_\C x})\,\bigr\}$.
\end{theorem}

The reverse inclusion ``$\supset$'' of this theorem was proved in
\cite[Theorem 5.6]{guillemin-sternberg;geometric-quantization} under
the assumption that $0$ is a regular value of $\Phi$, in \cite[Theorem
8.10]{kirwan;cohomology-quotients-symplectic} and in \cite[Section
2]{ness;stratification} under the assumption that $\omega$ is the
restriction to $M$ of the Fubini-Study symplectic form and $\Phi$ the
standard moment map, and in \cite[Proposition
2.4]{sjamaar;holomorphic} in the general case.

A second important case of Proposition \ref{proposition;length} is
when $\lie{k}=\b$, the Borel subalgebra of $\g_\C$.  Let
$\pr\colon\b\to\g$ be the restriction to $\b$ of the real linear
projection $\Im\colon\g_\C\to\g$.  Let $R$ be the root system of
$(G,T)$ and $R_+$ the set of positive roots.  Then
$$
\b=\t_\C\oplus\n=\t\oplus\sq\,\t\oplus\bigoplus_{\alpha\in
R_+}\g_\alpha.
$$
Hence $\pr$ sends $\t$ to $0$ and maps $\sq\,\t$ bijectively onto
$\t$.  Since $\bar{\g}_\alpha=\g_{-\alpha}$, $\pr$ maps $\n$
bijectively onto the $\Ad$-invariant complement
$$\h=\bigoplus_{\alpha\in R_+}(\g_\alpha\oplus\g_{-\alpha})\cap\g$$
of $\t$.  In particular, $\pr$ is onto.

The space of global sections is a $B$-module, so if $s\in\Gamma(M,L)$
is an eigensection for $\b$ with character $\chi\colon\b\to\C$, then
$\chi$ exponentiates to a character of $B$.  The two restriction maps
$$
\Hom(B,\C^\times)\longto\Hom(T_\C,\C^\times)\longto
\Hom(T,\U(1))=\Lambda
$$
are isomorphisms, because $B=T_\C N$ and $[B,B]=N$ is unipotent.
Therefore
$$
\chi(\xi)=2\pi\sq\,\lambda\bigl(\xi_1+\sq\,\xi_2\bigr)
=2\pi\sq\,\lambda(\xi_1)-2\pi\,\lambda(\xi_2)
$$
for a unique $\lambda\in\Lambda$.  (Here we decompose $\xi\in\b$ as
$\xi_1+\sq\,\xi_2+\xi_3$ with $\xi_1$, $\xi_2\in\t$ and $\xi_3\in\n$.)
This implies that $\Re\chi(\xi)=-2\pi\lambda(\xi_2)$.  Extending
$\lambda$ to an element of $\g^*$ by setting $\lambda=0$ on $\h$, and
recalling that $\pr\xi_3\in\h$, we obtain
$\Re\chi(\xi)=-2\pi\lambda(\pr\xi)$.  Thus Proposition
\ref{proposition;length} takes the form
$$
\ca{L}(\xi_M)\norm{s}^2
=4\pi\bigl(-\lambda(\pr\xi)+\phi(\pr\xi)\bigr)\norm{s}^2
$$
for all $\xi\in\b$ and all $N$-invariant holomorphic sections $s$ of
weight $\lambda$.  Similarly, for $s\in S^N_{r\lambda,r}$ we have
\begin{equation}\label{equation;semistable-unipotent}
\ca{L}(\xi_M)\norm{s}^2=4\pi
r\bigl(-\lambda(\pr\xi)+\phi(\pr\xi)\bigr)\norm{s}^2.
\end{equation}

\begin{proof}[Proof of Theorem {\rm\ref{theorem;maximum}}]
Let $x$ be a point in $X$ where the function $\norm{i^*s}$ attains its
maximum.  Then $\ca{L}(\xi_M)\norm{s}^2=0$ for all $\xi\in\b$, since
$X$ is $B$-invariant.  Moreover, $s(x)\ne0$ because $i^*s\ne0$.  Hence
$\phi(\pr\xi)(x)=\lambda(\pr\xi)$ by
\eqref{equation;semistable-unipotent}.  Since $\pr\colon\b\to\g$ is
surjective, this implies $\Phi(x)=\lambda$.
\end{proof}

\begin{remark}\label{remark;chamber}
Although for large $r$ the restriction map
$\Gamma(M,L^r)\to\Gamma(X,L^r)$ is surjective, the induced map on
$N$-invariants, $\Gamma(M,L^r)^N\to\Gamma(X,L^r)^N$, may not be.  (See
Remark \ref{remarks}\eqref{item;ideal}.)  Nevertheless, it follows
from Remark \ref{remark;norm}\eqref{item;subvariety} that Theorem
\ref{theorem;maximum} is also valid for sections $s\in\Gamma(X,L^r)^N$
that do not extend to global $N$-invariants.  Hence, if
$r\lambda\in\Lambda$ is any weight of $T_\C$ occurring in
$\Gamma(X,L^r)^N$, whether dominant or not, then $\lambda\in\Phi(X)$.
This proves the following assertion: let $A'(X)\supset A(X)$ be the
$\Lambda\times\N$-graded algebra $(S/I(X))^N$ and put
$$
\ca{C}'(X)=\bigl\{\,\lambda\in\Lambda\otimes\Q\bigm|
\text{$A'(X)_{r\lambda,r}\ne0$ for some $r>0$}\,\bigr\}.
$$
Then $\ca{C}'(X)$ is convex over $\Q$ and is a subset of
$\Phi(X)\cap\t^*$.
\end{remark}

Theorems \ref{theorem;maximum} and \ref{theorem;semistable} enable us
to establish one half of Theorem \ref{theorem;borel}.

\begin{theorem}\label{theorem;rational}
$\ca{C}(X)$ is equal to the set of rational points in $\Delta(X)$.
\end{theorem}

\begin{proof}
The inclusion $\ca{C}(X)\subset\Delta(X)\cap(\Lambda\otimes\Q)$
follows from Theorem \ref{theorem;maximum}.  Now let
$\lambda\in\Delta(X)$ be a rational point.  Then $\lambda=\Phi(x)$ for
some $x\in X$, and $r_0\lambda\in\Lambda_+$ for some positive integer
$r_0$.  Without loss of generality we may replace $r_0\lambda$ with
$\lambda$, $L$ with $L^{r_0}$, and $\Phi$ with $r_0\Phi$.  Then
$\Phi(x)=\lambda\in\Lambda_+$.  Let
$$
\p_\lambda=\b\oplus\bigoplus\bigl\{\,\g_{-\alpha}\bigm|
\text{$\alpha\in R_+$, $\lambda(\alpha\spcheck)=0$}\,\bigr\}
$$
be the parabolic subalgebra of $\g_\C$ associated with $\lambda$ and
$P_\lambda$ the parabolic subgroup of $G_\C$ generated by
$\exp\p_\lambda$.  Let $Y=G_\C/P_\lambda$ be the corresponding flag
variety and $L_{-\lambda}=(G_\C\times\C_{-\lambda})/P_\lambda$ the
Hermitian holomorphic line bundle on $Y$ determined by $-\lambda$.
Here $\C_{-\lambda}$ is the one-dimensional $P_\lambda$-module with
character $-\lambda$.  The curvature form $\omega_{-\lambda}$ of
$L_{-\lambda}$ is K\"ahler and a moment map $\Phi_{-\lambda}$ for the
$G$-action is found, as in \eqref{equation;flag}, by composing the
maps
$$
Y=G_\C/P_\lambda\longto G/G_\lambda\longto\g^*.
$$
Here $G_\lambda=P_\lambda\cap G$ is the centralizer of $\lambda$
(considered as an element of $\g^*$), the first map is the inverse of
the diffeomorphism $G/G_\lambda\to Y$ induced by the inclusion $G\to
G_\C$, and the second map is defined by $gG_\lambda\mapsto-g\lambda$.
Let us consider the product $M'=M\times Y$ with the K\"ahler form
$\omega'=\omega+\omega_{-\lambda}$ as a Hamiltonian $G$-manifold with
the diagonal $G$-action.  The moment map on $M'$ is
$\Phi'=\Phi+\Phi_{-\lambda}$, so from $\Phi(x)=\lambda$ we see that
$\Phi'(x,1P_\lambda)=0$.  Now $\omega'$ is the curvature form of the
Hermitian holomorphic line bundle $L'=L\boxtimes L_{-\lambda}$ and so,
by Theorem \ref{theorem;semistable}, there exist $r>0$ and
$s'\in\Gamma(M',(L')^r)^G$ satisfying
\begin{equation}\label{equation;bigsection}
s'(x,1P_\lambda)\ne0.
\end{equation}
By the Borel-Weil theorem, $\Gamma(Y,L_{-\lambda})\cong V_\lambda^*$,
the dual of the irreducible $G$-module $V_\lambda$ with highest weight
$\lambda$.  Accordingly, the K\"unneth formula yields isomorphisms of
$T_\C$-modules
\begin{multline*}
\Gamma(M',(L')^r)^G\cong
\bigl(\Gamma(M,L^r)\otimes\Gamma(Y,L_{-r\lambda})\bigr)^G\cong
\bigl(\Gamma(M,L^r)\otimes V_{r\lambda}^*\bigr)^G\\
\cong
\Hom(V_{r\lambda},\Gamma(M,L^r))^G\cong\Gamma(M,L^r)_{r\lambda}^N
=S^N_{r\lambda,r}.
\end{multline*}
A concrete isomorphism $\rho\colon\Gamma(M',(L')^r)^G\to
S^N_{r\lambda,r}$ is given by restricting invariant sections of
$(L')^r$ to the first factor $M\cong M\times\{1P_\lambda\}$, as
indicated in the pullback diagram
$$
\xymatrix{L^r\ar[d]\ar[r]&
L^r\boxtimes\C_{-r\lambda}\ar[r]\ar[d]&(L')^r\ar[d]\\
M\ar@/^1pc/[u]^s\ar[r]&M\times\{1P_\lambda\}\ar[r]
&**[r]M'\ar@/_1pc/[u]_{s'}.}
$$
Explicitly, if $s'\in\Gamma(M',(L')^r)^G$, then $\rho(s')$ is the
section of $L^r$ defined by $\rho(s')(y)\otimes1=s'(y,1P_\lambda)$ for
all $y\in M$.  It is easy to check that $\rho(s')$ is $N$-invariant
and transforms according to $r\lambda$ under the $T_\C$-action.
Conversely, a section $s\in S^N_{r\lambda,r}$ can be extended in a
unique $G_\C$-equivariant way to a global section $\epsilon(s)$ of
$(L')^r$, given by the formula
$\epsilon(s)(y,gP_\lambda)=g(s(g^{-1}y)\otimes1)$ for $y\in M$ and
$g\in G_\C$.  It is easy to verify that $\epsilon(s)$ is well-defined
and that the extension map $\epsilon\colon
S^N_{r\lambda,r}\to\Gamma(M',(L')^r)^G$ is the inverse of the
restriction map $\rho$.

Taking $s'$ as in \eqref{equation;bigsection} and putting $s=\rho(s')$
we find $s\in S^N_{r\lambda,r}$ satisfying $s(x)\ne0$.  Thus
$\lambda\in\ca{C}(X)$.  This proves
$\Delta(X)\cap(\Lambda\otimes\Q)\subset\ca{C}(X)$.
\end{proof}

%%%%%%%%%%%%%%%%%%%%%%%%%%%%%%%%%%%%%%%%%%%%%%%%%%%%%%%%%%%%%%%%%%%%%%%%
\section{Denseness of the rational points}\label{section;dense}
%%%%%%%%%%%%%%%%%%%%%%%%%%%%%%%%%%%%%%%%%%%%%%%%%%%%%%%%%%%%%%%%%%%%%%%%

In this section we will finish the proof of Theorem
\ref{theorem;borel} by establishing the following result.

\begin{theorem}\label{theorem;dense}
$\ca{C}(X)$ is dense in $\Delta(X)$.
\end{theorem}

The proof is a variation on the proof of Theorem
\ref{theorem;rational}, namely a semistability argument involving a
product $M\times Y_n$.  However, for our present purposes $Y_n$ will
not be a rational coadjoint orbit, but a member of a sequence
$(Y_n)_{n\ge1}$ of rational K\"ahler $G$-manifolds ``converging'' to
an irrational coadjoint orbit.

This auxiliary sequence of spaces is manufactured as follows.  Let
$\lambda$ be any point in the closed Weyl chamber $\t^*_+$.  The
chamber is a polyhedral cone and hence is a disjoint union of
(relatively open) faces.  Let $\sigma$ denote the unique face
containing $\lambda$, let $\t^*_\sigma$ be the linear span of $\sigma$
and let $\Lambda_\sigma$ be the lattice $\Lambda\cap\t^*_\sigma$.  The
centralizer $G_\sigma=G_\lambda$ of $\sigma$ is a connected subgroup
of $G$ containing the maximal torus $T$.  The root system of
$(G_\sigma,T)$ is the set $R_\sigma$ consisting of all roots
$\alpha\in R$ such that $\lambda(\alpha\spcheck)=0$, and so the
complexification of the Lie algebra $\g_\sigma$ has root-space
decomposition
$$
(\g_\sigma)_\C=\t_\C\oplus\bigoplus_{\alpha\in R_\sigma}\g_\alpha.
$$
From this we see that $\t^*_\sigma$ is the annihilator in
$\g_\sigma^*$ of the ideal $[\g_\sigma,\g_\sigma]$ of $\g_\sigma$, and
therefore is canonically isomorphic to the dual of the Lie algebra
$\t_\sigma$ of the torus $T_\sigma=G_\sigma/[G_\sigma,G_\sigma]$.
This implies that the canonical map
$$
\Hom(T_\sigma,\U(1))\overset\cong{\longto}\Hom(G_\sigma,\U(1))\longinj
\Hom(T,\U(1))=\Lambda
$$
sends the weight lattice of $T_\sigma$ isomorphically onto the
subgroup $\Lambda_\sigma$ of $\Lambda$.  We will use this natural
isomorphism to identify $\Hom(T_\sigma,\U(1))$ with $\Lambda_\sigma$.
Now let $(\Delta_n)_{n\ge1}$ be any sequence of polytopes in $\t^*$
subject to the following requirements:
\begin{enumerate}
\item\label{item;polytope-face}
$\Delta_n\subset\sigma$ for
all $n$;
\item\label{item;decreasing}
$\Delta_m\subset\Delta_n$ for $m\ge n$;
\item\label{item;limit}
$\bigcap_{n=1}^\infty\Delta_n=\{\lambda\}$;
\item\label{item;rational}
for each $n$ there exists a positive integer $d_n$ such that the
dilated polytope $d_n\Delta_n$ is a regular lattice polytope in
$\t^*_\sigma$ with respect to the lattice $\Lambda_\sigma$.
\end{enumerate}
(Recall that for any free abelian group $F$ a \emph{regular lattice
polytope} is a subset $P$ of the vector space $F\otimes_\Z\R$ which is
the convex hull of a finite subset of $F$ and each vertex $v$ of which
has the property that the cone over $P$ with apex $v$ is spanned by a
basis of $F$.)  Let $(Y^\circ_n,\omega^\circ_n)$ be the symplectic
toric $T_\sigma$-manifold associated with $\Delta_n$, as constructed
by Delzant \cite{delzant;hamiltoniens-periodiques}.  This is a compact
K\"ahler manifold equipped with a K\"ahlerian $T_\sigma$-action and a
$T_\sigma$-moment map $\Phi^\circ_n$, whose image is exactly
$\Delta_n$.  The form $d_n\omega^\circ_n$ is integral and so there
exists a unique $T$-equivariant Hermitian holomorphic line bundle
$L^\circ_n$ over $Y^\circ_n$ whose curvature form is equal to
$d_n\omega^\circ_n$.  Since $T_\sigma$ is just the abelianization of
$G_\sigma$, we can regard $Y^\circ_n$ as a Hamiltonian
$G_\sigma$-manifold and $L^\circ_n$ as a $G_\sigma$-equivariant line
bundle.  We can boost $Y^\circ_n$ to a Hamiltonian $G$-manifold $Y_n$
by the familiar process of \emph{symplectic induction}: put
$Y_n=(T^*G\times Y^\circ_n)\quot G_\sigma$, where $G_\sigma$ acts on
$T^*G$ by right multiplication and by the given action on $Y^\circ_n$,
and where the double slash indicates symplectic reduction at the zero
level.  The left $G$-action on $T^*G$ descends to a $G$-action on
$Y_n$ and the $G_\sigma$-moment map $\Phi^\circ_n$ induces a
$G$-equivariant map $\Phi_n\colon Y_n\to\g^*$, which is a moment map
for the $G$-action with respect to the reduced symplectic form
$\omega_n$.  The line bundle $L^\circ_n$ induces a Hermitian line
bundle $L_n$ on $Y_n$ with curvature form $d_n\omega_n$ and, finally,
the holomorphic structures on $Y^\circ_n$ and $L^\circ_n$ induce
$G$-invariant holomorphic structures on $Y_n$, resp.\ $L_n$.  (In
fact, there are isomorphisms of complex $G_\C$-manifolds
$$
Y_n\cong G_\C\times^{P_\lambda}Y^\circ_n\qquad\text{and}\qquad
L_n\cong G_\C\times^{P_\lambda}L^\circ_n,
$$
where the parabolic subgroup $P_\lambda$ acts on $Y^\circ_n$ and
$L^\circ_n$ via the Levi factor $(G_\lambda)_\C$.)  By construction,
the $G$-manifold $Y_n$ has moment polytope
$\Delta(Y_n)=\Phi_n(Y_n)\cap\t^*_+=\Delta_n$.

\begin{proof}[Proof of Theorem {\rm\ref{theorem;dense}}]
Let $\lambda\in\Delta(X)$.  We need to show that there exists a
sequence $(\lambda_n)_{n\ge1}$ in $\ca{C}(X)$ converging to $\lambda$.
We choose a sequence of polytopes $(\Delta_n)_{n\ge1}$ satisfying
requirements \eqref{item;polytope-face}--\eqref{item;rational} above
and we let $(Y_n)_{n\ge1}$ be the corresponding sequence of
Hamiltonian $G$-manifolds.  Let $Y_n^-$ denote the symplectic manifold
opposite to $Y_n$ and let $M'_n$ be the $G$-manifold $M\times Y_n^-$,
furnished with the K\"ahler form $\omega'_n=\omega-\omega_n$, the
moment map $\Phi'_n=\Phi-\Phi_n$ and the line bundle
$L'_n=L^{d_n}\boxtimes L_n^*$.  Note that the equivariant curvature
form of $L'_n$ is equal to $d_n(\omega'_n+\Phi'_n)$.  Choose any $x$
in $X$ such that $\Phi(x)=\lambda$.  Since
$\lambda\in\Delta_n=\Delta(Y_n)$, there exists $y_n\in Y_n$ such that
$\Phi'_n(x,y_n)=0$.  Therefore Theorem \ref{theorem;semistable}
(applied to the symplectic manifold $(M'_n,d_n\omega'_n)$) guarantees
the existence of a positive integer $r_n$ and a section
$s'_n\in\Gamma\bigl(M'_n,(L'_n)^{r_n}\bigr)^G$ such that
$s'_n(x,y_n)\ne0$.  Let $p_1\colon M'_n\to M$ and $p_2\colon M'_n\to
Y_n$ denote the Cartesian projections.  The natural maps
\begin{multline*}
\Gamma\bigl(M'_n,(L'_n)^{r_n}\bigr)\overset\cong{\longto}
\Hom\bigl(p_2^*L_n^{r_n},p_1^*L^{d_nr_n}\bigr)\\
\overset\cong{\longto}
\Hom\bigl(\Gamma(Y_n,L_n^{r_n}),\Gamma(M,L^{d_nr_n})\bigr)
\end{multline*}
are $G$-equivariant isomorphisms.  Hence, by Schur's lemma, there
exist a dominant weight $\mu$ and direct summands $E_\mu$ and $F_\mu$
of the $G$-modules $\Gamma(M,L^{d_nr_n})$, resp.\
$\Gamma(Y_n,L_n^{r_n})$, such that $E_\mu$ and $F_\mu$ are irreducible
of highest weight $\mu$ and such that the image of $s'_n$ in
$\Hom(F_\mu,E_\mu)^G$ does not vanish at the point $(x,y_n)$.
Therefore, if $s_E$ and $s_F$ are highest weight vectors in $E_\mu$,
resp.\ $F_\mu$, we have $s_E(x)\ne0$ and $s_F(y_n)\ne0$.  This shows
that the point $\mu/r_n$ is contained in both $\ca{C}(Y_n,L_n)$ and
$\ca{C}(X,L^{d_n})$.  In particular, by Theorem
\ref{theorem;rational},
$$
\mu/r_n\in\ca{C}(Y_n,L_n)\subset\Delta(Y_n,d_n\omega_n),
$$
and hence $\mu/r_n\in d_n\Delta_n$, since
$\Delta(Y_n,\omega_n)=\Delta_n$.  We conclude that the point
$\lambda_n=\mu/d_nr_n$ is contained in $\Delta_n\cap\ca{C}(X,L)$.
Requirements \eqref{item;decreasing}--\eqref{item;limit} now ensure
that $\lim_{n\to\infty}\lambda_n=\lambda$.
\end{proof}

A useful corollary of these results is the following extension of
Theorem \ref{theorem;semistable}.  Let us call a point $x\in M$
\emph{unipotently semistable} if there exist $r>0$ and
$s\in\Gamma(M,L^r)^N$ such that $s(x)\ne0$, in other words, if the set
$\ca{C}(\overline{Bx})$ is nonempty.  By Theorem \ref{theorem;dense},
$\Delta(\overline{Bx})$ is the closure of $\ca{C}(\overline{Bx})$, so
unipotent semistability can be characterized thus:

\begin{theorem}\label{theorem;unipotent-semistable}
A point $x$ is unipotently semistable if and only if
$\Delta(\overline{Bx})$ is nonempty.
\end{theorem}

%%%%%%%%%%%%%%%%%%%%%%%%%%%%%%%%%%%%%%%%%%%%%%%%%%%%%%%%%%%%%%%%%%%%%%%%
\section{``Schubert'' varieties}\label{section;schubert}
%%%%%%%%%%%%%%%%%%%%%%%%%%%%%%%%%%%%%%%%%%%%%%%%%%%%%%%%%%%%%%%%%%%%%%%%

In this section we will describe some subvarieties of $M$ to which
Theorem \ref{theorem;borel} applies.  Let $M^T_1$, $M^T_2$,\dots,
$M^T_k$ be the connected components of the $T$-fixed point set $M^T$.
For $j=1$, $2$,\dots, $k$, let $\alpha_{1j}$, $\alpha_{2j}$,\dots,
$\alpha_{c(j),j}$ be the nonzero weights of the isotropy
representation of the maximal torus $T$ on the tangent space at any
point of $M^T_j$, where $c(j)$ is the complex codimension of $M^T_j$.
Let $\ca{P}$ be the complement of the rational hyperplane arrangement
defined by this collection of weights:
$$
\ca{P}=\{\,\xi\in\t\mid\text{$\alpha_{ij}(\xi)\ne0$ for all $j=1$,
$2$,\dots, $k$ and $i=1$, $2$,\dots, $c(j)$}\,\}.
$$
The connected components of $\ca{P}$ are the \emph{action chambers} of
$\tau$.  We fix, once and for all, an action chamber $\ca{P}_+$ which
intersects the Weyl chamber $\t_+$.  Pick an element $\xi\in\ca{P}_+$
and put $\eta=\sq\,\xi\in\sq\,\t$.  For each $y\in M$ the limit
$x=\lim_{t\to-\infty}\exp(t\eta)y$ exists and is fixed under $\xi$.
Since $\xi\in\ca{P}$, this implies that $x\in M^T$.  The
\emph{unstable manifold} of $M^T_j$ is defined as the set
$$
W_j=\bigl\{\,y\in M\bigm|\lim_{t\to-\infty}\exp(t\eta)y\in
M^T_j\,\bigr\}.
$$
Generalizing results of Bia{\l}ynicki-Birula
\cite{bialynicki*;theorems-actions-algebraic-groups,
bialynicki*;properties-decompositions}, Carrell and Sommese
\cite{carrell-sommese;actions;scand;1978,
carrell-sommese;topological-aspects} showed that $W_j$ is a complex
submanifold of $M$, whose closure $X_j=\overline{W}_j$ is an analytic
subvariety in which $W_j$ is Zariski open.  Moreover, the map
$\pi_j\colon W_j\to M^T_j$ which sends $y$ to
$\lim_{t\to-\infty}\exp(t\eta)y$ is a holomorphic fibration, the
fibres of which are complex affine spaces.  Let us call $X_j$ a
\emph{generalized Schubert variety} of $M$.  The varieties $W_j$ and
$X_j$ and the fibration $\pi_j$ are defined by making a choice of an
element $\xi$ of the action chamber; however, one can show that they
don't depend on $\xi$ but only on the action chamber.  (See e.g.\
\cite[Section~3.6]%
{guillemin-sjamaar;convexity-properties-hamiltonian}.)  In particular,
we may, and will, assume that $\xi\in\ca{P}_+\cap\Int\t_+$.  Then the
parabolic subgroup associated with $\xi$ is equal to $B$, which
implies the following result.

\begin{proposition}\label{proposition;unstable}
Let $1\le j\le k$.  Then
\begin{enumerate}
\item
the unstable manifold $W_j$ and the generalized
Schubert variety $X_j$ are $B$-invariant;
\item
the map $\pi_j$ is $B$-invariant;
\item
for every $x\in M^T_j$ the fibre $W_{j,x}=\pi_j^{-1}(x)$ and its
closure $X_{j,x}=\overline{W}_{j,x}$ are $B$-invariant.
\end{enumerate}
\end{proposition}

\begin{proof}
It follows immediately from the definition that the submanifold $W_j$
and the map $\pi_j$ are $T_\C$-invariant.  Let
$$
\n=\bigoplus_{\alpha\in R_+}\g_\alpha
$$
be the decomposition of the Lie algebra of $N$ into weight spaces with
respect to the $T$-action, $R_+$ being the set of positive roots of
the root system $R$ of $(G,T)$.  Let $\zeta\in\n$ and write
$\zeta=\sum_{\alpha\in R_+}\zeta_\alpha$ with
$\zeta_\alpha\in\g_\alpha$.  Then $\Ad(\exp t\eta)\zeta=\sum_\alpha
e^{2\pi\alpha(\eta)t}\zeta_{\alpha}$.  Since $\xi\in\Int\t_+$, we have
$\alpha(\eta)>0$ and hence
$$
\zeta(t)=\Ad(\exp t\eta)\zeta\longto0\qquad\text{and}\qquad
\exp\zeta(t)\longto1
$$
as $t\to-\infty$.  In particular, if $y\in W_j$ and $n=\exp\zeta\in
N$,
$$
\exp(t\eta)ny=\exp(t\eta)n\exp(-t\eta)\exp(t\eta)y
=\exp\zeta(t)\exp(t\eta)y\longto\pi_j(y)
$$
as $t\to-\infty$.  This shows that $W_j$ and $\pi_j$ are
$N$-invariant.  Since $B=T_\C N$, both are $B$-invariant as well.  It
follows immediately that $X_j$, $W_{j,x}$ and $X_{j,x}$ are
$B$-invariant.
\end{proof}

This shows that one can apply Theorem \ref{theorem;borel} to the
generalized Schubert varieties $X_j$ and their subvarieties $X_{j,x}$.
As might be expected, for most $x$ in the fixed-point component
$M^T_j$ the polytopes $\Delta(X_j)$ and $\Delta(X_{j,x})$ are the
same.

\begin{corollary}\label{corollary;fibre-image}
Let $1\le j\le k$.  The set of all $x\in M^T_j$ such that
$\Delta(X_j)=\Delta(X_{j,x})$ is a nonempty Zariski open subset of
$M^T_j$.
\end{corollary}

\begin{proof}
Define
$$
U=\{\,y\in W_j\mid\Delta(\overline{By})=\Delta(X_j)\,\}
\quad\text{and}\quad C=\{\,x\in
M^T_j\mid\Delta(X_{j,x})=\Delta(X_j)\,\}.
$$
We claim that 
\begin{equation}\label{equation;image}
\pi_j(U)=C.
\end{equation}
Indeed, let $y\in U$ and put $x=\pi_j(y)$.  Then $\overline{By}$ is
contained in $X_{j,x}$ since $X_{j,x}$ is $B$-invariant, and therefore
$$
\Delta(X_j)=\Delta(\overline{By})\subset\Delta(X_{j,x})
\subset\Delta(X_j),
$$
i.e.\ $x\in C$.  This shows that $\pi_j(U)$ is a subset of $C$.
Conversely, let $x\in C$.  By Corollary \ref{corollary;orbit}, there
exists $y\in W_{j,x}$ such that
$\Delta(\overline{By})=\Delta(X_{j,x})$.  Hence
$\Delta(\overline{By})=\Delta(X_j)$, i.e.\ $y\in U$.  This shows that
$C$ is contained in $\pi_j(U)$, which proves \eqref{equation;image}.
By Corollary \ref{corollary;orbit}, $U$ is nonempty and Zariski open
in $W_j$ and therefore, by a theorem of Chevalley, its image
$\pi_j(U)=C$ is constructible in $M^T_j$.  Furthermore, $C$ is open in
the classical topology, since $\pi_j$ is a submersion.  Hence $C$ is
Zariski open.
\end{proof}

However, frequently the polytopes $\Delta(X_j)$ are not very exciting.
More interesting polytopes turn up when we restrict the $G$-action on
$M$ to a subgroup.  Particularly interesting is the case of the flag
variety $M=G_\C/B$.  Here the $W_j$'s are just the $B$-orbits in $M$,
so the $X_j$'s are the classical Schubert varieties.  Let us now view
$M$ as a $G_1$-space with respect to a closed connected subgroup $G_1$
of $G$.  Let us suppose $B$ is chosen in such a way that it contains a
Borel subgroup $B_1$ of $(G_1)_\C$, so that the $X_j$'s are
$B_1$-invariant.  If $G_1$ is the maximal torus $T$ of $G$, one gets
from Theorem \ref{theorem;borel} Atiyah's Theorem
\ref{theorem;schubert}, and if one takes $G_1$ to be an arbitrary
closed subgroup of $G$ one gets a nonabelian version of this theorem.

%%%%%%%%%%%%%%%%%%%%%%%%%%%%%%%%%%%%%%%%%%%%%%%%%%%%%%%%%%%%%%%%%%%%%%%%
%%%%%%%%%%%%%%%%%%%%%%%%%%%%%%%%%%%%%%%%%%%%%%%%%%%%%%%%%%%%%%%%%%%%%%%%

\bibliographystyle{amsplain}

\providecommand{\bysame}{\leavevmode\hbox to3em{\hrulefill}\thinspace}
\providecommand{\MR}{\relax\ifhmode\unskip\space\fi MR }
\providecommand{\MRhref}[2]{%
  \href{http://www.ams.org/mathscinet-getitem?mr=#1}{#2}
}
\providecommand{\href}[2]{#2}

%%%%%%%%%%%%%%%%%%%%%%%%%%%%%%%%%%%%%%%%%%%%%%%%%%%%%%%%%%%%%%%%%%%%%%%%
%%%%%%%%%%%%%%%%%%%%%%%%%%%%%%%%%%%%%%%%%%%%%%%%%%%%%%%%%%%%%%%%%%%%%%%%

\end{document}